\documentclass{amsart}
\usepackage{amssymb, latexsym}
\usepackage{amscd}

\newtheorem{theorem}{Theorem}
\newtheorem{lemma}{Lemma}
\newtheorem{remark}{Remark}
\newtheorem{definition}{Definition}

\newtheorem{corollary}{Corollary}

\begin{document}

\title{On Local Tameness of Certain Graphs of Groups}
\author{Rita Gitik}
\email{ritagtk@umich.edu}
\address{Department of Mathematics \\ University of Michigan \\ Ann Arbor, MI, 48109}  

\date{\today}

\begin{abstract} 
Let $G$ be the fundamental group of a finite graph of groups with Noetherian edges and locally tame vertices.
We prove that $G$ is locally tame. It follows that if a finitely presented group $H$ has a non-trivial $JSJ$-decomposition over the class of its $VPC(k)$ subgroups for $k=1$ or $k=2$, and all the vertex groups in the decomposition are flexible, then $H$ is locally tame.
\end{abstract}

\maketitle


 JSJ-decomposition, covering space, fundamental group.

\section{Introduction}

Let $H$ be a subgroup of a group $G$ given by the  presentation $G= \langle X|R \rangle $. 
Let $K$ be the standard presentation $2$-complex of $G$, i.e. $K$ has one 
vertex, $K$ has an edge, which is a loop, for every generator $x \in X$, and $K$ has a $2$-cell
for every relator $r \in R$. The Cayley complex of $G$, denoted by $Cayley_2(G)$, is  
the universal cover of $K$. Denote by $Cayley_2(G,H)$ the cover of $K$ corresponding to
a subgroup $H$ of $G$.

\begin{definition} cf. \cite{Gi1} and \cite{Mi1}.

A finitely generated subgroup $H$ of a finitely presented group  
$G$ is tame in $G$ if for any finite subcomplex
$C$ of $Cayley_2(G,H)$ and for any component $C_0$ of $Cayley_2(G,H)-C$ the group
$\pi_1(C_0)$ is finitely generated.
\end{definition}

A manifold $M$ is called a missing boundary manifold  if it can be 
embedded in a compact manifold $\bar M$  such that $\bar M - M$
is a closed subset of the boundary of $\bar M$. Simon conjectured in \cite{Sim} that if $M_0$ is a compact orientable irreducible $3$-manifold, and $M$ is the cover of $M_0$
corresponding to  a finitely generated subgroup of $\pi_1(M_0)$, then
$M$ is a missing boundary manifold.  Perelman's solution of Thurston's Geometrization Conjecture in 2003 implies that Simon's conjecture holds for all compact orientable irreducible $3$-manifolds, cf. \cite{B-B-B-M-P} and \cite{K-L}. 

Tucker proved in \cite{Tuc}  
that  a non-compact orientable irreducible $3$-manifold $M$ is a missing 
boundary manifold  if and only if the trivial subgroup is tame in
the fundamental group of $M$.

It is not known if there exists a finitely generated subgroup $H$ of a finitely presented group $G$ 
such that $H$ is not tame in $G$.

Tameness of a subgroup is connected to other properties which have been studied for a long time.

It is shown in \cite{M-T} that if the trivial subgroup is tame in $G$ then $\pi_1^\infty (G)$ (the fundamental group at infinity of $G$) is pro-finitely generated. It is shown in \cite{Mi1} that if a finitely generated subgroup $H$ is tame in $G$ then $\pi_1^\infty (G,H)$ is pro-finitely generated.

It is shown in \cite{M-T} that if the trivial subgroup is tame in $G$ then $G$ is QSF(Quasi-Simply-Filtrated).

\bigskip
 
The following definition was given in \cite{Gi2}.

\begin{definition}
A group $G$ is locally tame if all finitely generated subgroups of $G$
are tame in $G$.
\end{definition}

Recall that a group is called Noetherian or slender if all its subgroups are finitely generated. A group is polycyclic if it is Noetherian and solvable. For $n \ge 0$ a group $G$ is VPC(n), (virtually polycyclic of length $n$) if it has $n+1$ subgroups, $G_0, \cdots, G_n$ such that $G_{i+1}$ is a normal subgroup of $G_i$ for $0 \le i \le n-1$, the quotient groups $G_i/G_{i+1}$ are isomorphic to $\mathbf{Z}$ for $0 \le i \le n-1$, $G_n$ is the trivial subgroup, and $G_0$ has finite index in $G$.

Note that $VPC(0)$ groups are finite, $VPC(1)$ groups are finite extensions of $\mathbf{Z}$, and 
$VPC(2)$ groups are finite extensions of an extension of $\mathbf{Z}$ by   $\mathbf{Z}$. There are only two non-isomorphic extensions of $\mathbf{Z}$ by   $\mathbf{Z}$, namely the fundamental group of a torus and the fundamental group of a Klein bottle.

It is unknown whether all finitely presented Noetherian groups are virtually polycyclic (question 11.38 from the Kourovka Notebook, \cite{K-M}), however there exist finitely generated Noetherian groups that are not virtually polycyclic, for example the Tarski monster.

The main results of this paper is the following theorem.

\begin{theorem}
Let $G$ be a finitely presented group which is the fundamental group of a finite graph of groups with Noetherian edge groups.  If all the vertex groups of $G$ are locally
tame then $G$ is locally tame.
\end{theorem}

Recall that a subgroup $H$ is elliptic in a graph of groups $G$ if $H$ is contained in a conjugate of a vertex group. A vertex group $K$ of a $JSJ$-decomposition of $G$ which fails to be elliptic in some other $JSJ$-decomposition of $G$ is called flexible, cf. \cite{G-L}.

Theorem 1 implies the following interesting result.

\begin{lemma}
If a finitely presented group $G$ has a non-trivial $JSJ$-decomposition
over the class of its $VPC(k)$ subgroups for $k=1$ or $k=2$, and all the vertex groups in the decomposition are flexible, then $G$ is locally tame.
\end{lemma}

\begin{corollary}
Let $G$ be the fundamental group of a finite graph of groups which has all the vertex groups homeomorphic to $\mathbf{Z^n} \times ($surface group$)$ and all the edge groups homeomorphic to 
$\mathbf{Z^{n+1}}$. Then $G$ is locally tame.
\end{corollary}

\begin{remark} 
Let $G$ be a finitely presented group which has  a $JSJ$-decomposition over the class of its $VPC(n+1)$ subgroups. Let $K$  be a flexible  vertex group of this decomposition. Then $K$ is either $VPC(n+1)$ or $K$ has a finite index subgroup $L$ such that $L$ has a normal $VPC(n)$ subgroup $N$ with $L/N$  the fundamental group of a surface. Furthermore, if $L/N$ is the fundamental group of a closed surface, then $K=G$.
\end{remark}

\textbf{Conjecture.}
If a finitely presented group $G$ has a non-trivial $JSJ$-decomposition
over the class of its $VPC(n+1)$ subgroups for $n \ge 0$, and all the vertex groups in the decomposition are flexible, then $G$ is locally tame.

\section{Proof of Theorem 1}

We need the following notation.

Let $X^* = \{x,x^{-1} |x \in X \}$. For  $x \in X$ define $(x^{-1})^{-1} =x$.
 
Let $G$ be a group generated by a set $X$ and let $H$ be a subgroup of $G$.
Let $\{Hg \}$ be the set of right cosets  of $H$ in $G$.
 
The coset graph of $G$  with respect to $H$, denoted $Cayley(G,H)$,  is the oriented graph 
whose vertices are the cosets $\{Hg \}$, the set of edges is $ \{Hg \} \times X^*$, 
and an edge $ (Hg, x)$ begins at the vertex  $Hg$ and ends at the 
vertex $Hgx$. Denote the Cayley graph of $G$ by $Cayley(G)$. 
Note that $Cayley(G,H)$ is the quotient of $Cayley(G)$ by left 
multiplication by $H$. Also note that  the $1$-skeleton of $Cayley_2(G)$ is
$Cayley(G)$, and the $1$-skeleton of $Cayley_2(G,H)$ is $Cayley(G,H)$.

Let $G$ be generated by a disjoint union of sets $X_i, 1 \le i \le n$. 
We call a connected subcomplex $C$ of $Cayley(G,H)$ an $X_i$-component, if all edges of $C$
have the form $(Hg, x)$ with $x \in X_i^*$.

\textbf{Proof of Theorem 1.}

Let $G$ be a finite graph of groups with vertex groups $V_i, 1 \le i \le n$ and edge groups $E_j, 1 \le j \le m$.
As $G$ is finitely presented and all the edge groups are Noetherian, hence finitely generated, it follows that all the vertex groups are finitely presented.
Let the vertex group $V_i$ be generated by a finite set $X_i$ such that the sets $X_i$ and $X_k$ are disjoint for $i \neq k$.

Consider a finitely generated subgroup  $H$ of $G$. Note that $H$ is the fundamental group of a (possibly infinite) graph of groups which has the vertex groups isomorphic to subgroups 
of conjugates of $V_i$ and the edge groups isomorphic to subgroups of conjugates of $E_j$, \cite{S-W}. 

As the edge groups of $G$ are Noetherian, the edge groups of $H$ are also Noetherian and the vertex groups of $H$ are finitely generated.

Note that all maximal $X_i$-components of $Cayley_2(G,H)$ have fundamental groups which are subgroups of conjugates of $V_i$, hence the  maximal $X_i$-components of $Cayley_2(G,H)$ are homeomorphic 
to $Cayley_2(V_i, U_i)$, with $U_i$  a finitely generated subgroup of $V_i$. 

As $H$ is finitely generated, there exists a finite connected subcomplex
$(K, H \cdot 1)$ of $Cayley_2(G,H)$ such that the inclusion map of $(K, H \cdot 1)$ in $Cayley_2(G,H)$ induces an isomorphism of $\pi_1(K, H \cdot 1)$ with $\pi_1(Cayley_2(G,H), H \cdot 1)=H$.

Let $C$ be a compact subcomplex of $Cayley_2(G,H)$. Note that $C$ has non-empty intersection
with only finitely many maximal $X_i$-components of $Cayley_2(G,H)$.  The complex $K$ can be enlarged to contain $C$. It can be enlarged more, so it consists of finitely many maximal $X_i$-components of $Cayley_2(G,H)$ which have non-trivial intersection with $C$ and the $2$-cells with boundaries in the union of those $X_i$-components. By construction, $K-C$ has a finite number of connected components.

As the vertex groups $V_i$ are locally tame, the fundamental group of each component of the complement of $C$ in any maximal $X_i$-component is finitely generated,  hence the fundamental group of each component of $K-C$ is finitely generated,  

Note that $(Cayley_2(G,H) - C)=(Cayley_2(G,H)-K)\cup (K-C)$. Let $W$ be a connected component of the closure of $Cayley_2(G,H) -K$. Then $W \cap K$ is connected and $\pi_1(W \cap K)$ is isomorphic to $\pi_1(W)$ because $K$ carries the fundamental group of $Cayley_2(G,H)$. So for each component $K_i$ of $K-C$ which intersects $W$ non-trivially, 
$\pi_1(K_i \cap W)=\pi_1(W)$. Let $W^i$ be the (possibly infinite) union of all components of $Cayley_2(G,H)-K$ which have non-trivial intersection with $K_i$. Then $\pi_1(W^i \cup K_i)=\pi_1(K_i)$ which is finitely generated. Hence the fundamental group of each component of $Cayley_2(G,H) - C$ is finitely generated, proving Theorem 1.

\section{Proof of Lemma 1}

\begin{remark}
The following result was proved in \cite{Gi2}.
Let $K_0$ be a finite index subgroup of a finitely presented group $K$. A finitely generated 
subgroup $H$ of $K$ is tame in $K$ 
if and only if $H \cap K_0$ is tame in $K_0$.

It follows that virtually locally tame groups are locally tame. 
\end{remark}

\begin{remark}
Note that the fundamental group of a surface is locally tame.

It is shown in \cite{Gi2} that finitely generated free groups are locally tame. Indeed, for any free group $F$ and its finitely generated subgroup $H$ the complex $Cayley_2(F,H)$
is one-dimensional. When $H$ is finitely generated, $Cayley_2(F,H)$ is homotopic to
a wedge of finitely many circles. It follows that the fundamental group of a non-closed surface is tame.

It is shown in \cite{Gi2} that finitely generated abelian groups are locally tame, hence the fundamental group of a torus is locally tame.

Note that the fundamental group of a closed orientable surface of genus greater than one can be written as a double of a free group over a cyclic subgroup. Hence Theorem 1 implies that fundamental groups of closed orientable surfaces of genus greater than one are locally tame. 

As closed orientable surfaces are double covers of non-orientable closed surfaces of the same genus,
Remark 2 implies that the fundamental groups of non-orientable closed surfaces are locally tame.
\end{remark}

\textbf{Proof of Lemma 1.}

Consider, first, the case when a finitely presented group $G$ has a non-trivial $JSJ$-decomposition
over the class of its $VPC(1)$ subgroups and all the vertex groups in the decomposition are flexible. Note that $VPC(1)$ groups are Noetherian.

The flexible vertex groups in such $JSJ$-decomposition are either  $VPC(1)$ or virtually(fundamental group of surfaces), cf. \cite{S-S} and \cite{G-L}. Furthermore, if a vertex group $M$ in that decomposition is virtually(the fundamental group of a closed surface), then $G=M$.

Hence Remark 2 and Remark 3 imply that the group $G$ satisfies the conditions of Theorem 1, therefore it is locally tame.
 
Next, consider the case when a finitely presented group $G$ has a non-trivial $JSJ$-decomposition
over the class of its $VPC(2)$ subgroups and all the vertex groups in the decomposition are flexible. Note that $VPC(2)$ subgroups are Noetherian.

The flexible vertex groups in such $JSJ$-decomposition are either $VPC(2)$ or virtually-(cyclic-by-a surface group), cf. \cite{S-S} and \cite{G-L}. Furthermore, if a flexible vertex group $K$ in that decomposition is virtually-(cyclic-by-a closed surface group), then $G=K$.

If a group $L$ is (cyclic-by-a surface group) then there exists   a surface $M$ and a normal cyclic subgroup $N$ of $L$ such that the following  sequence is exact.
$$1 \rightarrow N \rightarrow L \rightarrow \pi_1(M) \rightarrow 1 $$
and $L$ is the fundamental group of a bundle $X$ over $M$ with fiber $\mathbf{S^1}$.

If $H$ is a finitely generated subgroup of $L$ then either $H \cap N = \{ 1\}$ or $H \cap N$ is isomorphic to $\mathbf{Z}$. Let $K$ be the image of $H$ in $\pi_1(M)$. Note that $K$ is finitely generated. Let $M_K$ be the cover of $M$ with fundamental group $K$. Then $H$ is the fundamental group of a bundle $X_H$ over $M_K$ with fiber either $\mathbf{S^1}$ if $H \cap N = Z$ or fiber $\mathbf{R}$ if
 $H \cap N = \{ 1\}$. As $K$ is finitely generated, $M_K$ is a missing boundary surface. It follows that, in either case, $X_H$ is a missing boundary $3$-manifold, so $L$ is locally tame.

If  a group $L$ is $VPC(2)$ then it is virtually either the fundamental group of a torus or the fundamental group of a Klein bottle, hence Remark 3 implies that $L$ is locally tame.

Therefore, Remark 2 implies that the group $G$ satisfies the conditions of Theorem 1, so it is locally tame.

\section{Acknowledgement}
The author would like to thank Mike Mihalik and Peter Scott for helpful discussions.

\end{document}